\documentclass[onefignum,onetabnum]{siamart171218}

\usepackage{lipsum}
\usepackage{amsfonts}
\usepackage{graphicx}
\usepackage{epstopdf}

\usepackage{multirow}

\usepackage{tikz}
\usetikzlibrary{matrix}
\usetikzlibrary{shapes,patterns,snakes,arrows}
\usetikzlibrary{positioning,fit,calc}
\usetikzlibrary{graphs}

\usepackage{algorithmic}
\ifpdf
  \DeclareGraphicsExtensions{.eps,.pdf,.png,.jpg}
\else
  \DeclareGraphicsExtensions{.eps}
\fi

\newsiamremark{remark}{Remark}
\newsiamremark{hypothesis}{Hypothesis}
\crefname{hypothesis}{Hypothesis}{Hypotheses}
\newsiamthm{claim}{Claim}

\def\bb{\begin{bmatrix}}
\def\eb{\end{bmatrix}}

\newcommand{\TheTitle}{Machine learning and domain decomposition methods - a survey}

\headers{Machine learning and domain decomposition methods}{A. Klawonn, M. Lanser, and J. Weber}

\title{{\TheTitle}}

\author{Axel Klawonn\footnotemark[1]\ \footnotemark[2] \and Martin Lanser\footnotemark[1]\ \footnotemark[2] \and Janine Weber\footnotemark[1]\ \footnotemark[2]}

\usepackage{amsopn}

\usepackage{ifthen}%
\newcommand{\mycomment}[1]{%
\ifthenelse{\isodd{\value{page}}}{%
\normalmarginpar%
\marginpar{\tiny {#1}}%
}{%
\reversemarginpar%
\marginpar{\tiny {#1}}%
}}%

\usepackage{cite}
\usepackage[english]{babel}
\usepackage{hyperref} 
\usepackage{amsmath}
\usepackage{amssymb}
\usepackage{amsfonts}
\usepackage[T1]{fontenc}
\usepackage[utf8]{inputenc}

\usepackage{makeidx}  
\usepackage{enumerate}
\usepackage{tikz}
\usetikzlibrary{matrix}
\usetikzlibrary{shapes,patterns,snakes,arrows,decorations.markings,arrows.meta}
\usetikzlibrary{positioning,fit,calc}
\usetikzlibrary{graphs}
\usepackage{ifthen}
\usepackage{nicefrac}
\usepackage{relsize}
\usepackage{bbm}

\usepackage{listings}  
\usepackage{booktabs}

\definecolor{darkspringgreen}{rgb}{0.09, 0.45, 0.27}
\definecolor{darkpastelgreen}{rgb}{0.01, 0.75, 0.24}
\definecolor{bostonuniversityred}{rgb}{0.8, 0.0, 0.0}
\definecolor{ferrarired}{rgb}{1.0, 0.11, 0.0}
\definecolor{coralred}{rgb}{1.0, 0.25, 0.25}
\definecolor{pastelred}{rgb}{1.0, 0.41, 0.38}
\definecolor{amber}{rgb}{1.0, 0.75, 0.0}
\definecolor{darkblue}{RGB}{0,0,139}
\definecolor{darkred}{RGB}{105,0,0}

\definecolor{mgreen}{HTML}{017101}
\definecolor{dblue}{HTML}{004D80}
\definecolor{lblue}{HTML}{00AC8E}

\def\x{{\bf x}}

\newcommand{\VGDSW}[1]{V_{\rm{GDSW}}}

\newcommand{\takeout}[1]{ }

\tikzset{%
  every neuron/.style={
    circle,
    draw,
    minimum size=0.6cm
  },
  every input neuron/.style={
    circle,
    draw,
    minimum size=0.6cm,
    fill=green!50
  },
  every output neuron/.style={
    circle,
    draw,
    minimum size=0.6cm,
    fill=orange!30
  },
  every hidden neuron/.style={
    circle,
    draw,
    minimum size=0.6cm,
    fill=blue!40
  },
  neuron missing/.style={
    draw=none, 
    scale=1.5,
    text height=0.3cm,
    execute at begin node=\color{black}$\vdots$
  },
}

\begin{document}

\maketitle

\renewcommand{\thefootnote}{\fnsymbol{footnote}}

\footnotetext[1]{Department of Mathematics and Computer Science, University of Cologne, Weyertal 86-90, 50931 K\"oln, Germany, \email{\{axel.klawonn, martin.lanser, janine.weber\}@uni-koeln.de}, url: \url{http://www.numerik.uni-koeln.de}} 
\footnotetext[2]{Center for Data and Simulation Science, University of Cologne, Germany, url: \url{http://www.cds.uni-koeln.de}}

\begin{abstract}
Hybrid algorithms, which combine black-box machine learning methods with experience from traditional numerical methods and domain expertise from diverse application areas, are progressively gaining importance in scientific machine learning and various industrial domains, especially in computational science and engineering. 
In the present survey, several promising avenues of research will be examined which focus on the combination of machine learning (ML) and domain decomposition methods (DDMs). The aim of this survey is to provide an overview of existing work within this field and to structure it into domain decomposition for machine learning and machine learning-enhanced domain decomposition, including: domain decomposition for classical machine learning, domain decomposition to accelerate the training of physics-aware neural networks, machine learning to enhance the convergence properties or computational efficiency of DDMs, and machine learning as a discretization method in a DDM for the solution of PDEs.
In each of these fields, we summarize existing work and key advances within a common framework and, finally, disuss ongoing challenges and opportunities for future research. 
\end{abstract}



\section{Introduction}
\label{sec:intro}
Domain decomposition methods (DDMs) are divide-and-con\-quer strategies which decompose a given problem into a number of smaller subproblems. This strategy often leads to easier parallelizable algorithms where the subproblems can be solved on different processors (cpu or gpu). Sometimes, an additional problem is needed for parallel scalability, the so-called global problem, or to gather and connect local information obtained from the solution of the local subproblems. DDMs  can be useful in different cases, for example, if the original problem is too large to be solved in the available memory, the computing time is too long on a single processor, or the decomposed problem has other preferable properties, for example, being better conditioned, yielding more accurate results, or boosting generalization properties and mitigating the spectral bias, respectively, of machine learning models. 
DDMs have a long and successful history in constructing parallel scalable preconditioners for linear and nonlinear systems obtained from discretized partial differential equations. Recently, there has been an increased interest in using DDMs in combination with machine learning algorithms, especially in the field of scientific machine learning (SciML).

\begin{figure}
\centering
\includegraphics[width=0.9\textwidth]{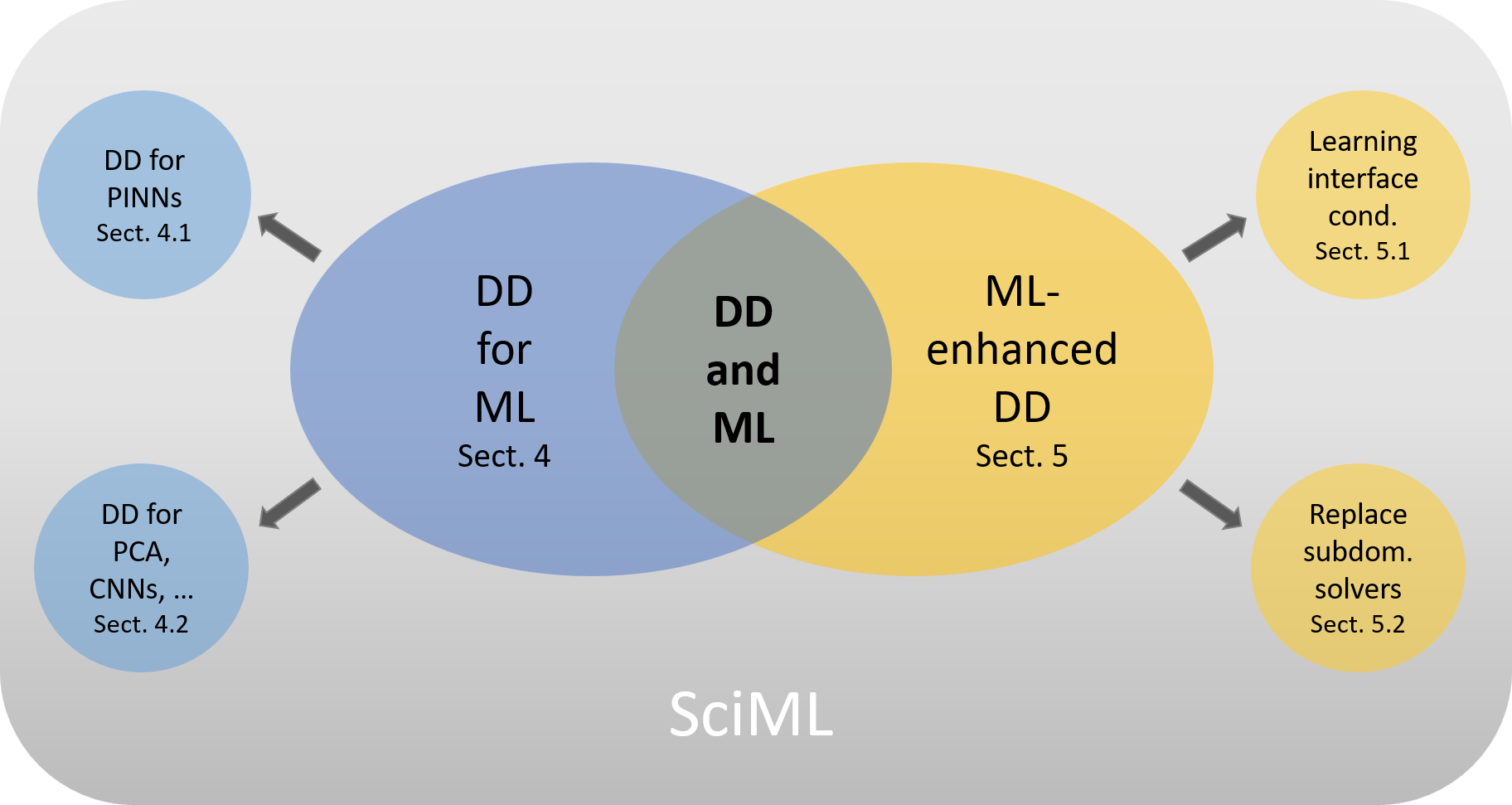}
\caption{\label{fig:struct} Schematic representation of the structure of this review paper. All reviewed work on the combination of machine learning and domain decomposition methods are clustered into the two main groups domain decomposition (DD) for machine learning (ML),~\cref{sec:dd_in_ml}, and machine learning-enhanced domain decomposition,~\cref{sec:ml_in_dd}.}
\end{figure}

In recent years, the rapidly growing research area of SciML has drawn increasing interest and attention in various fields of applications, as, for example, in computational fluid mechanics~\cite{grimm2023learning,eichinger2022surrogate}, optics and electromagnetic applications~\cite{chen2020physics,kovacs2022conditional}, and earthquake detections~\cite{smith2020eikonet}. Originally, the term SciML was introduced and made publicly known by the report~\cite{USDpt:SciML} prepared for the U.S. Department of Energy. 
 The core idea of SciML is to combine algorithms from supervised or unsupervised machine learning with other domain-specific methods or knowledge to develop new, hybrid methods. In particular, the authors of~\cite{USDpt:SciML} have identified six priority research directions in SciML which are, among others, domain-awareness, interpretability, and robustness of the considered method. With respect to practical applications, this means that a key factor of SciML methods is to connect ``black-box'' approaches such as deep learning with well-understood and theoretical approaches from science and expert knowledge. 
 
One framework that fits very well to the priority research directions of SciML are physics-informed neural networks (PINNs)~\cite{RPK:2019:MID}. PINNs are a specific machine learning technique usually used to solve forward or inverse problems involving Partial Differential Equations (PDEs). Similar to other deep learning models, PINNs approximate PDE solutions by training a neural network such that a given loss function is minimized. However, the characteristic feature of PINNs is that knowledge with respect to the underlying physics behind the data is integrated into the loss function as well as, possibly, initial or boundary
conditions of the computational domain's boundary. For a detailed overview of existing research and applications around PINNs as well as remaining challenges and open research questions, we refer to~\cite{cuomo2022scientific}.

One drawback of PINNs and also other neural network architectures, in general, is the large computational complexity of the related optimization problem for large problem domains and multi-scale problems as well as the resulting long training times.
Hence, a wide range of attempts have been made to distribute and parallelize the training of different network architectures; see~\cite{ben-nun_demystifying_2020,peteiro2013survey,verbraeken2020survey} for a comprehensive overview. 
In general, the wide range of techniques used in parallel and distributed deep learning can roughly be categorized into \textit{data parallelism}, that is, partitioning the input samples, \textit{model parallelism}, that is, partitioning the network structure, and \textit{pipelining}, that is, partitioning by layer; see also~\cite[Sect. 6]{ben-nun_demystifying_2020}. 
From an abstract perspective, model parallelism can be and often is interpreted as a form of DDMs~\cite{toselli,QuarteroniValli2008}.  

In this paper, our aim is to review existing work on combining DDMs with different machine learning models. Although data parallelism can be also interpreted as a form of domain decomposition such that the global data space is decomposed into smaller subspaces operating only on parts of the data, in the following, we focus on model parallelism and pipelining when using domain decomposition in machine learning algorithms. 
 The remainder of the paper is organized as follows. In~\cref{sec:dd}, we provide a general description of DDMs for the solution of PDEs and introduce some necessary notation. Subsequently, in~\cref{sec:ml_pinns}, we briefly summarize the main idea and central algorithmic formulae of multilayer perceptrons, PINNs, and the Deep Ritz method. Note that both,~\cref{sec:dd} and~\cref{sec:ml_pinns}, can be skipped by the experienced reader who is already familiar with DDMs and PINNs and who can directly proceed with~\cref{sec:dd_in_ml}. 
In the present paper, we categorize all cited work into two main classes:
\begin{enumerate}[i)]
	\item Methods using domain decomposition within machine learning models (\cref{sec:dd_in_ml})
	\item Machine learning-enhanced DDMs (\cref{sec:ml_in_dd}). 
\end{enumerate}

Additionally, we classify all methods in group i) into two subgroups. First, we consider all approaches using domain decomposition for the acceleration and parallel training of PINNs; see~\cref{sec:dd_in_pinns}. Second, in~\cref{sec:dd_in_others}, we consider approaches using domain decomposition in machine learning algorithms other than PINNs, in particular, classical supervised and unsupervised algorithms. 
With respect to work within class ii), one can also subdivide the approaches into two subclasses; see also~\cite{HKLW:2020:GammReview}. In this case, the first class consists of methods where machine learning is used to improve the convergence properties or the computational efficiency within classical DDMs. This is typically done by learning or approximating optimal interface conditions or by learning other optimal parameters; see~\cref{sec:ml_in_dd_interface}.
In the second class, different types of neural networks are used as discretization methods and replace classical local subdomain or coarse solvers based on finite elements or finite differences; see~\cref{sec:ml_in_dd_solvers}.
For a schematic overview of the described categorization and the structure of this paper, see~\cref{fig:struct}. 

Note that for some approaches which combine domain decomposition with machine learning the categorization is not uniquely defined. For example, D3M~\cite{Li:2019:D3M} and DeepDDM (Deep Domain Decomposition Method)~\cite{Li:2020:DeepDDM} can be interpreted as replacing subdomain solvers within a DDM by neural networks but, at the same time, also as a DDM within neural network training. Hence, we have made the following classifications to the best of our knowledge and with regard to the subject focus of the concrete considered work.

\section{Domain decomposition methods for the solution of PDEs}
\label{sec:dd}

\begin{figure}[t]
\centering
\scalebox{0.7}{
\begin{tikzpicture}[x=0.75cm, y=0.75cm, line width=1pt]
	
	
	\foreach \j in {0,...,2} {
		\foreach \i in {0,...,2} {
			\ifthenelse{\i=\j \OR \(\i=2 \AND \j=0\) \OR \(\i=0 \AND \j=2\)}{\draw[line width=1pt,fill=darkpastelgreen,opacity=0.5] (3*\i,3*\j) -- (3*\i+3,3*\j) -- (3*\i+3,3*\j+3) -- (3*\i,3*\j+3) -- cycle;}{\draw[line width=1pt,fill=amber,opacity=0.5] (3*\i,3*\j) -- (3*\i+3,3*\j) -- (3*\i+3,3*\j+3) -- (3*\i,3*\j+3) -- cycle;}
		
		}
	}
	

	\draw[] (1.5,1.5) node[fill=white] {$\Omega_{1}$};
	\draw[] (4.5,1.5) node[fill=white] {$\Omega_{2}$};
	\draw[] (7.5,1.5) node[fill=white] {$\Omega_{3}$};
	
	\draw[] (1.5,4.5) node[fill=white] {$\Omega_{4}$};
	\draw[] (4.5,4.5) node[fill=white] {$\Omega_{5}$};
	\draw[] (7.5,4.5) node[fill=white] {$\Omega_{6}$};
	
	\draw[] (1.5,7.5) node[fill=white] {$\Omega_{7}$};
	\draw[] (4.5,7.5) node[fill=white] {$\Omega_{8}$};
	\draw[] (7.5,7.5) node[fill=white] {$\Omega_{9}$};

	\draw[->,>=stealth,blue](3,4.5) -- node[above] {$\delta$} (2,4.5);
	\draw[->,>=stealth,blue](6,4.5) -- node[above] {$\delta$} (7,4.5);
	\draw[->,>=stealth,blue](4.5,3) -- node[left] {$\delta$} (4.5,2);
	\draw[->,>=stealth,blue](4.5,6) -- node[left] {$\delta$} (4.5,7);

	\draw[] (5.5,5.5) node[blue,fill=white] {$\Omega_{5}'$};	
	\draw[line width=2pt,fill=blue!50!gray,opacity=0.5] (2,2) rectangle (7,7);
	\draw[line width=2pt,blue] (2,2) rectangle (7,7);

	\foreach \j in {0,...,2} {
		\foreach \i in {0,...,2} {
			\draw[line width=1pt,bostonuniversityred] (11+2*\i,2*\j+1.5) -- (11+2*\i+2,2*\j+1.5) -- (11+2*\i+2,2*\j+3.5) -- (11+2*\i,2*\j+3.5) -- cycle;
		}
	}
	
	\draw[line width=1pt] (11,1.5) -- (17,1.5) -- (17,7.5) -- (11,7.5) -- cycle;
	
	\draw[] (11+1,2.5) node[fill=white] {$\Omega_{1}$};
	\draw[] (11+3,2.5) node[fill=white] {$\Omega_{2}$};
	\draw[] (11+5,2.5) node[fill=white] {$\Omega_{3}$};
	
	\draw[] (11+1,4.5) node[fill=white] {$\Omega_{4}$};
	\draw[] (11+3,4.5) node[fill=white] {$\Omega_{5}$};
	\draw[] (11+5,4.5) node[fill=white] {$\Omega_{6}$};
	
	\draw[] (11+1,6.5) node[fill=white] {$\Omega_{7}$};
	\draw[] (11+3,6.5) node[fill=white] {$\Omega_{8}$};
	\draw[] (11+5,6.5) node[fill=white] {$\Omega_{9}$};

	\node[circle,fill=darkpastelgreen,inner sep=0pt,minimum size=6pt] (V1) at (13,3.5) {};
	\node[circle,fill=darkpastelgreen,inner sep=0pt,minimum size=6pt] (V2) at (15,3.5) {};
	\node[circle,fill=darkpastelgreen,inner sep=0pt,minimum size=6pt] (V3) at (13,5.5) {};
	\node[circle,fill=darkpastelgreen,inner sep=0pt,minimum size=6pt] (V4) at (15,5.5) {};
	
	\draw[] (17.5,3.5) node[fill=white] {\color{bostonuniversityred} $\Gamma$};
	
\end{tikzpicture}
}		
\caption{Decomposition of a domain $\Omega \subset \mathbb{R}^2$ into nine nonoverlapping subdomains $\Omega_i,\;i=1,...,9$. {\bf Left:} The overlapping domain decomposition is obtained by extending the nonoverlapping subdomains by the width $\delta$. The resulting overlapping subdomain $\Omega_5'$ is marked in blue. {\bf Right:} The interface $\Gamma$ is marked in red. }
\label{fig:dd}
\end{figure}
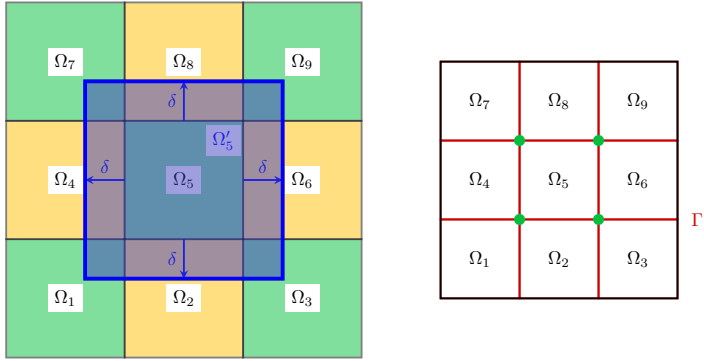

In this section, we briefly introduce the main principles of overlapping and nonoverlapping DDMs as well as some necessary notation.

DDMs~\cite{toselli,QuarteroniValli2008} are highly scalable, iterative methods for the solution of sparse, linear systems of equations as, for example, arising from the discretization of PDEs. DDMs have been shown to be numerically stable for many practically relevant  
problems and are designed for the application on parallel computers. Generally speaking, they rely on a divide-and-conquer strategy and decompose the original global problem into a number of smaller subproblems. 
Mathematically, this corresponds to a subdivision of the domain into a number of smaller subdomains, where the local problems can be solved in parallel on different processors of a parallel computer. Additionally, a small global problem is often needed for numerical and parallel scalability. Certain mathematical conditions at the interface of neighboring subdomains have to be satisfied such that the solution of the original problem is recovered.

For a generic description of DDMs, we consider a boundary value problem of the general form
\begin{equation}
\begin{array}{ll}
\mathcal{L}(u) &= f \quad {\rm in} \quad \Omega\\
\mathcal{B}(u) &= g \quad {\rm on} \quad \partial  \Omega
\end{array}
\label{eq:bvp}	
\end{equation}
on the domain $\Omega \subset \mathbb{R}^{d},\; d=2,3$, where $\mathcal{L}$  
is a linear, second-order, elliptic differential operator and $\mathcal{B}$ represents the boundary conditions.

In order to compute a numerical solution of boundary value problem~\eqref{eq:bvp}, we discretize the variational formulation of the given linear PDE with an appropriate numerical method, for example, conforming finite elements. 
This results in a, usually sparse, linear system of equations
\begin{equation}
K_g u_g = f_g. 
\label{eq:lin_system}	
\end{equation}
The basic idea of DDMs is then, instead of directly solving the completely assembled system of equations~\eqref{eq:lin_system}, to decompose the problem into smaller subproblems and thus, solving smaller systems of equations which are assembled on local parts of the domain exclusively.
 In the following, we denote these local parts of the domain as subdomains and we assume to have a decomposition of $\Omega$ into $N \in \mathbb{N}$ nonoverlapping subdomains $\Omega_i,\; i=1,...,N$, which fulfills $\overline{\Omega} = \bigcup_{i=1}^N \overline{\Omega}_i$; see also~\cref{fig:dd}. 
 Moreover, we define $\Gamma := \left(\bigcup_{i=1}^N \partial \Omega_i\right)\setminus \partial \Omega$ as the interface between neighboring subdomains, that is, degrees of freedom shared by at least two subdomains. 
 In order to define overlapping DDMs, we additionally introduce overlapping subdomains $\Omega'_i,\; i=1,...,N$, where the overlap $\delta$ is a measure for the amount of local information shared between neighboring subdomains; see also~\cref{fig:dd} (left) for a visualization. 

However, simply solving decoupled local subdomain systems independently from each other might result in discontinuities along the interface $\Gamma$ and will not result in the correct solution of the original and global system. Hence, in order to obtain the correct global solution which is continuous across the different subdomains, also communication between the different subproblems is necessary, that is, a global coupling between the different subdomains has to be ensured. With respect on how the global exchange of information across the original domain $\Omega$ is ensured, we distinguish between one- and two-level DDMs.  
One-level DDMs, where information is only exchanged between neighboring subdomains, result in a rate of convergence which deteriorates with a growing number of subdomains when used for the approximate solution of elliptic systems of PDEs; see, for example,~\cite{toselli}.
 A remedy is obtained using two-level DDMs by the setup and the solution of an additional, small global problem which needs to be solved in each iteration of the iterative solver. 
In particular, in two-level DDMs, the additional global problem ensures the fast global transport of information between the different subdomains as well as the scalability and robustness of the iterative method.
In the context of two-level DDMs, we refer to this global problem as \textit{coarse problem} and to the related solution space as \textit{coarse space}. 
Let us note that the idea of a coarse problem is also transferred to domain decomposition approaches for, e.g., training PINNs; cf.~\cref{sec:dd_in_pinns}.
 In classic DDMs, a coarse space is often defined by certain weighted averages or translations and rotations along the interface between neighboring subdomains. Additionally, also more sophisticated coarse spaces exist which are, for example, obtained by solving certain eigenvalue problems and which are problem-dependent or, in other words, adaptive. However, a detailed overview of existing classic and adaptive coarse spaces in DDMs would be beyond the scope of this review article.

\section{Multilayer perceptrons, PINNs, and Deep Ritz}
\label{sec:ml_pinns}

In this section, we briefly provide some preliminaries with respect to multilayer perceptrons, including a short outlook on convolutional neural networks and residual neural networks, see~\cref{sec:nn}, as well as PINNs and the Deep Ritz method; see~\cref{sec:pinns}. For more details on deep learning and more specialized network architectures, we refer to, for example,~\cite{Goodfellow:2016:DL,haykin2010neural,chollet2017deep}.
 Let us note that for experienced readers, who are already familiar with multilayer perceptrons, PINNs, and the Deep Ritz method, it might be more convenient to skip the following section and proceed directly with~\cref{sec:dd_in_ml}.

\subsection{Feedforward multilayer perceptrons}
\label{sec:nn}

\begin{figure}[t]
\centering
\scalebox{0.85}{
\begin{tikzpicture}[line width=1pt]
\def\xoffset{4.2}
\def\yoffset{0.4}

\foreach \m [count=\y] in {1,2,missing,j,missing,n} {
  \if\y3
  	\node [every neuron/.try, neuron \m/.try] (input-\m) at (\xoffset-0.4,\yoffset+3.85-0.8*\y) {};
  \else
  \if\y5
  	\node [every neuron/.try, neuron \m/.try] (input-\m) at (\xoffset-0.4,\yoffset+3.85-0.8*\y) {};
  \else
  	\node [every input neuron/.try, neuron \m/.try] (input-\m) at (\xoffset-0.4,\yoffset+3.85-0.8*\y) {};
  \fi
  \fi
  }

\foreach \m [count=\y] in {1,2,missing,K} {
  \if\y3
  	\node [every neuron/.try, neuron \m/.try ] (hidden1-\m) at (\xoffset+1.6,\yoffset+3.5-\y*1.0) {};
  \else
  	\node [every hidden neuron/.try, neuron \m/.try ] (hidden1-\m) at (\xoffset+1.6,\yoffset+3.5-\y*1.0) {\small{$h_\m^1$}};
  \fi
  }
  
\foreach \m [count=\y] in {1,2,missing,K}
  \node (hidden2-\m) at (\xoffset+3.2,\yoffset+3.5-\y*1.0) {};

\foreach \m [count=\y] in {1,2,3,4}
  \node at (\xoffset+3.7,\yoffset+3.5-\y*1.0) {\dots};
  
\foreach \m [count=\y] in {1,2,missing,K}
  \node (hidden3-\m) at (\xoffset+4.0,\yoffset+3.5-\y*1.0) {};

\foreach \m [count=\y] in {1,2,missing,K} {
  \if\y3
  	\node [every neuron/.try, neuron \m/.try ] (hidden4-\m) at (\xoffset+5.6,\yoffset+3.5-\y*1.0) {};
  \else
  	\node [every hidden neuron/.try, neuron \m/.try ] (hidden4-\m) at (\xoffset+5.6,\yoffset+3.5-\y*1.0) {\footnotesize{$h_\m^ N$}};
  \fi
  }  

\foreach \m [count=\y] in {1,missing,m} {
  \if\y2
  	\node [every neuron/.try, neuron \m/.try] (output-\m) at (\xoffset+7.2+0.4,\yoffset+2.4-0.75*\y) {};
  \else
  	\node [every output neuron/.try, neuron \m/.try] (output-\m) at (\xoffset+7.2+0.4,\yoffset+2.4-0.75*\y) {};
  \fi
  }

\foreach \i in {1,2,j,n}
  \foreach \j in {1,2,K}
    \draw [->] (input-\i) -- (hidden1-\j);
    
\foreach \i in {1,2,K}
  \foreach \j in {1,2,K}
    \draw [->] (hidden1-\i) -- (hidden2-\j);

\foreach \i in {1,2,K}
  \foreach \j in {1,2,K}
    \draw [->] (hidden3-\i) -- (hidden4-\j);    

\foreach \i in {1,2,K}
  \foreach \j in {1,m}
    \draw [->] (hidden4-\i) -- (output-\j);
    
\foreach \i in {1,2,j,n}
  \draw [<-] (input-\i) -- ++(-1,0)
    node [above, midway] {$i_{p,\i}$};  
    
\foreach \i in {1,m}
  \draw [->] (output-\i) -- ++(1,0)
    node [above, midway] {$o_{p,\i}$};      

\foreach \l [count=\x from 0] in {Input\\layer, Hidden\\layers, Output\\layer}
  \node [align=center, above] at (\xoffset-0.4+\x*4.0,4.0) {\l};

\end{tikzpicture}
}
\caption[Structure of a dense feedforward neural network with $n$ input nodes, $m$ output nodes, and $N$ hidden layers with $K$ neurons per layer.]{Structure of a dense feedforward neural network with $n$ input nodes (marked in green), $m$ output nodes (marked in orange), and $N$ hidden layers with $K$ neurons per layer (marked in blue). Figure taken from~\cite[Fig. 7.1]{Weber:2022:Diss}.}
\label{fig:fnn}
\end{figure}
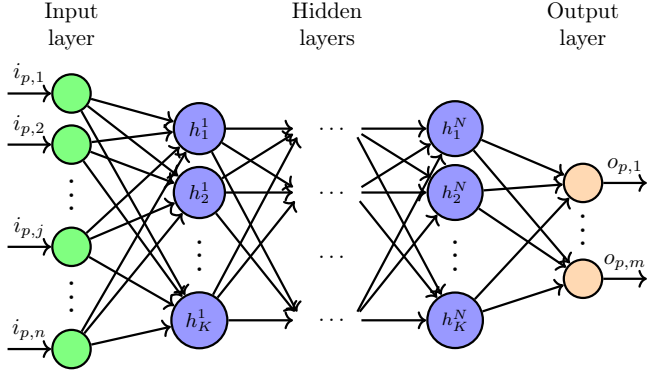

In general, neural networks or multilayer perceptrons, respectively, are used within supervised learning to approximate a functional relation between given input and output data. Usually, the overall aim is to minimize a loss function which measures the difference between the ground truth of the given data and the output, that is, estimated classification or regression values, of the neural network. 
The following description of multilayer perceptrons is loosely based on~\cite{HKLW:2018:ML_acc,Shalev-Shwartz:2014:UML,Weber:2022:Diss}. 
Generally, a feedforward neural network can be interpreted as a directed, weighted graph ${\cal G} = ({\cal V},{\cal E})$ with a set of nodes ${\cal V}$, a set of edges ${\cal E}$, and a weight function $w : {\cal E} \rightarrow \mathbb{R}$; see, for example, \cite[Chapt.~20.1]{Shalev-Shwartz:2014:UML}. An exemplary visualization of the graph of a dense feedforward neural network is presented in~\cref{fig:fnn}.
 A neural network is assumed to be organized in layers, that is, the set of nodes $\cal V$ can be represented as the union of nonempty, disjoint subsets ${\cal V}_i \subset {\cal V}, i= 0,\ldots, N+1$. For feedforward neural networks, these sets are defined such that for each edge $e \in \cal E$ there exists an $i \in \{0, \ldots, N\}$ with $e$ being an edge between a node in ${\cal V}_{i}$ and one in ${\cal V}_{i+1}$; see \cite[Chapt.~20.1]{Shalev-Shwartz:2014:UML}. 
 In general, different layers can perform different transformations on their input. A neural network usually starts with an input layer (marked in green in~\cref{fig:fnn}), where the features $\{ i_p\}_{p=1}^P$ of the external data are used as input data, and concludes with an output layer (marked in orange in~\cref{fig:fnn}). We denote the layers in between the input and the output layer as hidden layers and the nodes in a neural network as neurons. 

Mathematically, the relation between two consecutive layers is the 
conjunction of a linear mapping defined by the weight function $w$ and a nonlinear activation function $\alpha$. In particular, the nonlinearity of the activation function enables a neural network to approximate highly complex relations between the input and output data.

Generally, for a generic activation function $\alpha$, the output of the $k$-th layer of the neural network can be 
written as
\begin{align*}
	y = \alpha^k(x,W^k,{ b}^k),
\end{align*}
where $W^k { = (w^k_{ij})_{i,j}}$ and $b^k$ are the  weight matrix and the bias vector, 
respectively. Note that an entry $w^k_{ij}$ of $W^k$ corresponds to the value of the weight function $w$ associated with the corresponding edge between layer $ {\cal V}_{k-1}$ and $ {\cal V}_k$. Then, the application of a complete neural network 
with $N$ hidden layers to an input vector $i \in I$ is given by
\begin{equation}
\label{eq:nn_layers}
\begin{aligned}
	h^1 	& = \alpha^1(i,W^1,{ b^1}), \\
	h^{k+1} & = \alpha^{ k+1}(h^{ k},W^{k+1},b^{k+1}), \quad 1 { \leq} k < N, \\
	o 		& = { (W^{N+1})^T} h^{N} + { b^{N+1}}
\end{aligned}
\end{equation}
where  
$h^k$  
 is the output of the $k$-th hidden layer and $o \in O$ is the (final) output vector. The computation of the output vector $o$ is usually performed without an additional application of the activation function. 
  Note that for specific neural network architectures, the consecutive relation between two layers can be more complicated than written in~\eqref{eq:nn_layers}. In~\cref{sec:dd_in_others}, we consider also DDMs for the training of convolutional neural networks (CNNs) and residual neural networks (ResNets). In CNNs, convolutions, that is, kernel functions, and pooling functions are gradually applied to input data with a grid-like structure such that the information within the data is aggregated and condensed throughout the network layers. For more mathematical details on CNNs, we refer to, for example,~\cite{lecun1989generalization},~\cite[Chapt. 9]{Goodfellow:2016:DL}, and~\cite[Chapt. 5]{chollet2017deep}. For ResNets, skip connections that perform identity mappings are added between distant blocks of hidden layers of a multilayer perceptron or a CNN. Hence, the connected block of hidden layers approximates the residual function with respect to the input data which can help to mitigate the degradation problem of very deep neural networks; see~\cite{he2016deep} for further details.

\subsection{PINNs and Deep Ritz}
\label{sec:pinns}

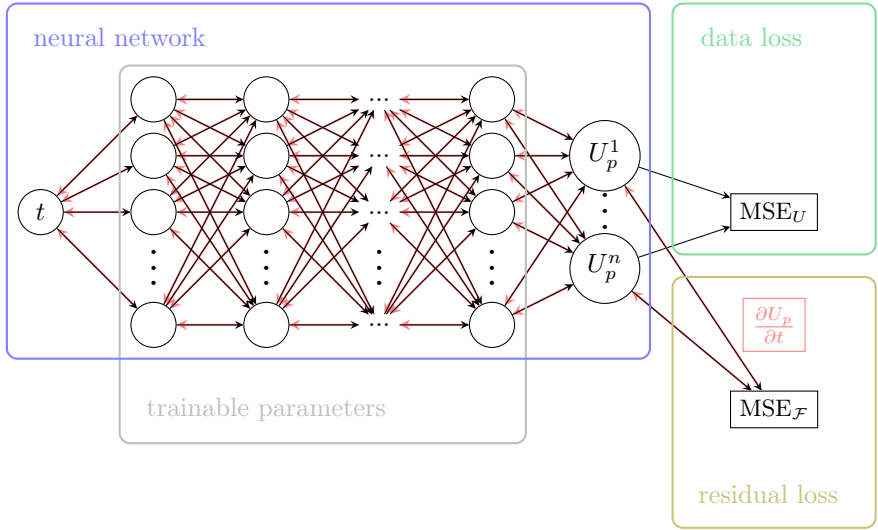
\begin{figure}[t]
\centering
\begin{tikzpicture}[x=1.5cm, y=1.5cm, >=stealth]

\node [every neuron/.try] (input-1) at (0,1.0) {$t$};

\foreach \m [count=\y] in {1,2,3,missing,4}
  \node [every neuron/.try, neuron \m/.try] (hiddenI-\m) at (1,2.5-0.5*\y) {};
  
\foreach \m [count=\y] in {1,2,3,missing,4}
  \node [every neuron/.try, neuron \m/.try] (hiddenII-\m) at (2,2.5-0.5*\y) {};
  
\foreach \m [count=\y] in {1,2,3,missing,4} {
\if\y4
	\node [every neuron/.try, neuron \m/.try] (empty-\m) at (3,2.5-0.5*\y) {};
\else
	\node (empty-\m) at (3,2.5-0.5*\y) {...};
\fi
}

\foreach \m [count=\y] in {1,2,3,missing,4}
  \node [every neuron/.try, neuron \m/.try] (hiddenIII-\m) at (4,2.5-0.5*\y) {};

\node [every neuron/.try, neuron 1/.try ] (output-1) at (5,1.5-0.5*0) {$U_p^1$};
\node [every neuron/.try, neuron missing/.try ] (output-missing) at (5,1.5-0.5*1) {};
\node [every neuron/.try, neuron 2/.try ] (output-2) at (5,1.5-0.5*2) {$U_p^{n}$};

\node[draw,rectangle] (mse-u) at (6.5,1.0) {\small $\text{MSE}_{U}$};
\node[draw,rectangle] (mse-f) at (6.5,-0.75) {\small $\text{MSE}_\mathcal{F}$};

\foreach \i in {1}
  \foreach \j in {1,...,4}
    \draw[red!50,thick,->] (hiddenI-\j) -- (input-\i);

\foreach \i in {1,...,4}
  \foreach \j in {1,...,4}
    \draw[red!50,thick,->] (hiddenII-\j) -- (hiddenI-\i);
    
\foreach \i in {1,...,4}
  \foreach \j in {1,...,4}
    \draw[red!50,thick,->] (empty-\j) -- (hiddenII-\i);
    
\foreach \i in {1,...,4}
  \foreach \j in {1,...,4}
    \draw[red!50,thick,->] (hiddenIII-\j) -- (empty-\i);
    
\foreach \i in {1,...,4}
  \foreach \j in {1,2}
    \draw [red!50,thick,->] (output-\j) -- (hiddenIII-\i);

\foreach \i in {1}
  \foreach \j in {1,...,4}
    \draw[->] (input-\i) -- (hiddenI-\j);

\foreach \i in {1,...,4}
  \foreach \j in {1,...,4}
    \draw[->] (hiddenI-\i) -- (hiddenII-\j);
    
\foreach \i in {1,...,4}
  \foreach \j in {1,...,4}
    \draw[->] (hiddenII-\i) -- (empty-\j);
    
\foreach \i in {1,...,4}
  \foreach \j in {1,...,4}
    \draw[->] (empty-\i) -- (hiddenIII-\j);
    
\foreach \i in {1,...,4}
  \foreach \j in {1,2}
    \draw[->] (hiddenIII-\i) -- (output-\j);

\draw[->] (output-1) -- (mse-u);

\draw[->] (output-2) -- (mse-u);

\draw[red!50,thick,->] (mse-f) -- (output-1);
\draw[->] (output-1) -- (mse-f);

\draw[red!50,thick,->] (mse-f) -- (output-2);
\draw[->] (output-2) -- (mse-f);

\node[draw,red!50,rectangle] at (6.5,0.0) {\color{red!50} $\frac{\partial U_p}{\partial t}$};

\draw[gray!50,rounded corners,thick] (0.7,-1.05) rectangle (4.3,2.3);
\node (train-par) at (2.0,-0.75) {\color{gray!50} trainable parameters};

\draw[blue!50,rounded corners,thick] (-0.3,-0.3) rectangle (5.4,2.85);
\node (train-par) at (0.7,2.55) {\color{blue!50} neural network}; 

\draw[darkpastelgreen!50,rounded corners,thick] (5.6,0.625) rectangle (7.4,2.85);
\node (train-par) at (6.3,2.55) {\color{darkpastelgreen!50} data loss};

\draw[olive!50,rounded corners,thick] (5.6,-1.8) rectangle (7.4,0.425);
\node (train-par) at (6.45,-1.5) {\color{olive!50} residual loss};

\end{tikzpicture}
\caption{Graphical representation of a PINN. The derivatives needed for the evaluation of the residual are computed using automatic differentiation in a backward propagation (red). Figure in modified form in~\cite[Fig. 3]{grimm2022Corsim}.}
\label{fig:pinn}
\end{figure}

The key concept of PINNs is to integrate domain specific knowledge in form of physical laws or domain expertise modeled by ordinary or partial differential equations into neural networks.  
In particular, this is done by enhancing the loss function with the residual error of a certain differential equation; see~\cite{RPK:2019:MID} for more details. Hence, the objective is to obtain solutions which do not only fit some given training data but also satisfy a given ordinary or partial differential equation in a least square sense. The following descriptions are loosely based on~\cite[Sect. 4]{grimm2022Corsim}.

 For the mathematical derivation of a physics-enhanced loss function, we consider again the general boundary value problem in~\eqref{eq:bvp} on the domain $\Omega \subset \mathbb{R}^{d},\; d=2,3$ with a linear, second-order, elliptic differential operator $\mathcal{L}$ and boundary conditions $\mathcal{B}$. 
We now aim to define a PINN which solves~\eqref{eq:bvp}, that is, we aim for a feedforward neural network $\mathcal{N}(\cdot,W,b)$ that fulfills~\eqref{eq:bvp} in a least square sense with weights $W$ and biases $b$. As input data of the neural network, we use a set of collocation points located inside the domain $\Omega$ as well as a set on the boundary $\partial \Omega$. In order to obtain a neural network $\mathcal{N}$ solving the boundary value problem~\eqref{eq:bvp}, we enhance the classic loss function by a point-wise error of the residual of the PDE.  Thus, the physics-informed loss function is defined as
\begin{equation}
\begin{array}{ll}
\mathcal{M}(W,b) &:= \mathcal{M}_\Omega(W,b) + \mathcal{M}_{\partial \Omega}(W,b)\\
\mathcal{M}_\Omega(W,b) &:= \frac{1}{N_f} \sum\limits_{i=1}^{N_f} |\mathcal{L}(\mathcal{N}(x_f^{i},W,b))-f(x_f^{i})|^2\\
\mathcal{M}_{\partial\Omega}(W,b) &:= \frac{1}{N_g} \sum\limits_{i=1}^{N_g} |\mathcal{B}(\mathcal{N}(x_g^{i},W,b))-g(x_g^{i})|^2,
\end{array}
\label{eq:loss}	
\end{equation}
with collocation points $x_f^{i},\; i=1,...,N_f,$ located in the domain $\Omega$ and collocation points $x_g^{i},\; i=1,...,N_g,$ located on the boundary $\partial \Omega$. 
During the training of the neural network, the derivatives occurring in the operators $\mathcal{L}$ and $\mathcal{B}$ are evaluated using the backward propagation algorithm and automatic differentiation~\cite{Baydin:2017:AD}.
Taking advantage of this, the training of the neural network then consists of solving the optimization problem
\begin{equation}
\{W^*,b^*\} := {\rm arg}\min\limits_{\{W,b\}} \mathcal{M}(W,b)
\label{eq:opt}	
\end{equation}
using a stochastic gradient approach based on mini-batches built from all collocation points. An exemplary PINN is shown in~\cref{fig:pinn}.
Note that the loss term $\mathcal{M}_\Omega(W,b)$ enforces the PINN to satisfy the condition $\mathcal{L}(u)=f$ in a least square sense, while $\mathcal{M}_{\partial \Omega}(W,b)$ enforces the boundary condition.
 
An approach that is very closely related to PINNs is the Deep Ritz method~\cite{Weinan:2018:DeepRitz} which uses dense neural networks (DNNs) for solving PDEs in variational form.  
The Deep Ritz method is based on the Ritz approach to formulate a given PDE as in~\eqref{eq:bvp} as an equivalent minimization problem which is then discretized and solved by a DNN in combination with a numerical integration method.
For a formal description of the Deep Ritz method, we assume that~\eqref{eq:bvp} is self-adjoint in the following. Then, solving the PDE is equivalent to the solution of the minimization problem
\begin{equation}
\min_{u} \mathcal{E}(u) \quad \text{s.t.} \quad \mathcal{B}(u) = g \quad {\rm on} \quad \partial  \Omega 
\label{eq:ritz}
\end{equation}
 where $\mathcal{E}(u)$ is an appropriate energy functional.
The basic idea of the Deep Ritz method is now to approximate the function $u$ in~\eqref{eq:ritz} by a DNN and to use a numerical quadrature rule where each integration point in $\Omega$ is used as a collocation point to approximate the minimizing functional $\mathcal{E}(u)$. 
 This yields the discrete loss function for the Deep Ritz method 
\begin{equation}
\begin{aligned}
\label{eq:loss_vp3}
\widetilde{\mathcal{M}}(W,b) &:= \widetilde{\mathcal{M}}_\Omega(W,b) + \widetilde{\mathcal{M}}_{\partial \Omega}(W,b) \quad \text{with}\\
\widetilde{\mathcal{M}}_{\Omega}(W,b)&  = \frac{1}{\widetilde{N}_f} \sum_{i=1}^{\widetilde{N}_f} h(\widetilde{x}_f^{i},W,b) \\
\widetilde{\mathcal{M}}_{\partial \Omega}(W,b) & = q\ \frac{1}{\widetilde{N}_g} \sum_{i=1}^{\widetilde{N}_g} \vert \mathcal{B}(\widetilde{\mathcal{N}}(\widetilde{x}_g^{i},W,b))-g(\widetilde{x}_g^{i})\vert^2,
\end{aligned}
\end{equation}
 where $h(\cdot)$ is the inner function of the energy functional $\mathcal{E}(u)$ approximated by the DNN $\widetilde{\mathcal{N}}$ and $q$ is a Lagrange multiplier.
 The derivatives of the neural network $\widetilde{\mathcal{N}}$ occuring in~\cref{eq:loss_vp3} are then evaluated using the backward propagation algorithm and automatic differentiation~\cite{Baydin:2017:AD}.  For a detailed derivation of the loss function~\eqref{eq:loss_vp3} for the  exemplary case of $\mathcal{L}(u) = \Delta u$, we refer to~\cite[Sect. 3.2]{HKLW:2020:GammReview}.

\section{Domain decomposition for machine learning}
\label{sec:dd_in_ml}

In this section, we review recent works where ideas from domain decomposition are used to accelerate and parallelize the training of different machine learning models, and, in particular, deep learning models.
 As described in~\cref{sec:intro}, the wide range of techniques used in parallel and distributed deep learning can roughly be categorized into \textit{data parallelism}, \textit{model parallelism}, and \textit{pipelining}; see also~\cite[Sect. 6]{ben-nun_demystifying_2020}. 
 Let us note again that, within this paper, we focus on model parallel training methods and pipelining rather than on data parallelism; see also~\cref{sec:intro}. In particular, the interpretation of model parallelism as a DDM is, from our point of view, more straightforward and also allows for analogous concepts in distributed ML corresponding to a coarse space in DDMs.

In the following, we will distinguish the different methods using DDMs in ML into two classes. First, we consider all approaches using domain decomposition for the acceleration and parallel training of very specific neural networks, that is, physics-informed neural networks, in~\cref{sec:dd_in_pinns}. For details on PINNs, we refer to~\cref{sec:pinns}. 
Second, in~\cref{sec:dd_in_others}, we provide an overview of approaches using domain decomposition in other ML algorithms than PINNs, in particular, classical supervised and unsupervised algorithms.

\subsection{Domain decomposition to accelerate the training of physics-aware neural networks}
\label{sec:dd_in_pinns}

\begin{table}
\scalebox{0.837}{
\begin{tabular}{l||l|l|l}
DD for PINNs & DD & PINN formulation & Interface conditions \\\hline
{\color{mgreen} cPINN~\cite{jagtap_conservative_2020} } & nonoverlapping & conservation law & conservative quantity \\
{\color{mgreen}  XPINN~\cite{karniadakis_extended_2020,hu_when_2022 } } & nonoverlapping & residuals & averaging \\ 
{\color{lblue} GatedPINN~\cite{stiller_large-scale_2020} } & adaptive & residuals & weighted average \\
{\color{lblue} APINN~\cite{hu2023augmented} } & adaptive & residuals & weighted average \\
{\color{mgreen} DPINN~\cite{dwivedi_distributed_2019} } & nonoverlapping & residuals & flux interface loss \\
{\color{mgreen} PECANN~\cite{basir_generalized_2023} } & nonoverlapping & residuals & generalized Robin-type \\
{\color{dblue}  Mosaic Flow~\cite{wang2022mosaic,feeney_breaking_2023} } & overlapping  & residuals & alterating Schwarz-like \\
{\color{dblue} PINN-PINN, PINN-FOM~\cite{snyder2023domain} } & overlapping  & residuals & loss term/strong enforc. \\
{\color{dblue} Iterative algorithms~\cite{yang_iterative_2023} } & overlapping  & residuals & weighted average \\
{\color{mgreen} $hp$-VPINNs~\cite{kharazmi_hp-vpinns_2021} } & nonoverlapping & var. residuals & via variational loss\\
{\color{mgreen} TgNN-wf~\cite{kharazmi_hp-vpinns_2021} } & nonoverlapping & var. residuals & via variational loss \\
{\color{dblue} PFNN-2~\cite{sheng_pfnn-2_2022} } & overlapping & var. residuals & via boundary network \\
{\color{dblue} FBPINN~\cite{moseley_finite_2023,dolean_finite_2023} } & overlapping & residuals & Schwarz-like \\
{\color{mgreen}  Nonlinear preconditioning~\cite{kopanicakova_enhancing_2023} } & nonoverlapping & residuals & additive/multiplicative \\
\end{tabular}
}
\caption{{\label{tab:DDinPINNs}} Summary of the current state of domain decomposition methods to accelerate the training of PINNs as reviewed in~\cref{sec:dd_in_pinns}. We categorize all papers according to the defined domain decomposition, the considered PINN formulation, and the concrete method of how to integrate interface constraints between neighboring subdomains.
We mark the papers according to the following color code: dark green: based on a nonoverlapping domain decomposition, blue: based on an overlapping domain decomposition, light green: an adaptive domain decomposition automatically trained by one or several neural networks.}
\end{table}

In general, integrating ideas from domain decomposition into PINNs usually aims at accelerating and parallelizing the training of PINNs, mainly by reducing the complexity of the PINN optimization problem for large problem domains and multi-scale problems. Roughly speaking, this means that instead of solving one large optimization problem, the idea is to use a divide-and-conquer strategy by solving, that is, training smaller local problems in parallel. This helps to parallelize the training and, potentially, can also help to reduce the complexity of the global optimization problem. 
Analogously to DDMs, the resulting related key question within these parallel training approaches is how to ensure that the different local problems communicate with each other such that the local solutions match across the interface. In the following, we review different parallel training methods for different types of PINNs. We will observe that some approaches deal with the communication between different local problems by integrating additional interface terms in the PINN loss function, whereas other approaches define overlapping subdomains and directly exchange information among the overlapping subdomain areas. In~\cref{tab:DDinPINNs}, a structured overview of all papers reviewed in the following can be found.

In~\cite{jagtap_conservative_2020}, conservative PINNs (cPINNs) are introduced that make use of a nonoverlapping 
domain decomposition approach for conservation laws. Here, the spatial computational domain is decomposed into nonoverlapping subdomains and for each subdomain, a separate local PINN, which also enforces certain interface conditions with its neighboring subdomains with respect to the solution and the local PDE residuals, is trained.  
In~\cite{jagtap_conservative_2020}, the authors present numerical results for different scalar nonlinear and systems of conservation laws.

The idea of cPINNs is further extended in~\cite{karniadakis_extended_2020} to generic PDEs and arbitrary space-time domains resulting in extended PINNs (XPINNs). Similar as in~\cite{jagtap_conservative_2020}, local sub-PINNs are trained for local subdomains of a domain independently from each other and in parallel. In particular, in XPINNs, the domain can be decomposed in both, space and time and no explicit assumptions of a specific PDE type are made. The approximated solutions of the local sub-PINNs are combined by enforcing continuity across the interface by averaging local solutions across neighboring subdomain interfaces as well as by enforcing a residual continuity condition. 
 In~\cite[Sect. 3.2]{de2022error}, the authors propose to additionally enforce continuity of first-order derivatives between different sub-PINNs to further reduce the training error close to the interface.

First steps regarding a direct comparison of PINNs and XPINNs with respect to generalization properties are presented in~\cite{hu_when_2022}. Specifically, the authors derive an a priori generalization bound based on the complexity of the considered PDE problem as well as an a posteriori generalization bound based on the matrix norm of the neural network after optimization and use this theory to derive observations with respect to the generalization properties of XPINNs. On the one hand, the decomposition of the global problem into smaller subproblems in XPINNs decreases the complexity of each subproblem and hence can boost generalization. On the other hand, each subproblem is trained with a smaller amount of training data which tends to increase overfitting of the local models. In~\cite{hu_when_2022}, different experimental results are presented which validate that the two factors in the derived generalization bounds lead to a tradeoff when comparing the performance of XPINNs to PINNs.

In~\cite{shukla_parallel_2021}, a unified parallel algorithm is presented for cPINNs and XPINNs which employs domain decomposition in space and space-time, respectively. The distributed framework in~\cite{shukla_parallel_2021} for the specific PINN extensions from~\cite{jagtap_conservative_2020} and~\cite{karniadakis_extended_2020} is constructed by a hybrid programming model using MPI~\cite{Sniretal1998,Groppetal1998}. Within the unified framework, the authors compare the parallel performance of cPINNs and XPINNs for different forward problems. The experimental results indicate that cPINNs are more effective regarding communication cost whereas XPINNs provide greater flexibility such that they can also handle space-time domain decompositions for any differential equations as well as deal with arbitrarily shaped subdomains.

The authors of~\cite{stiller_large-scale_2020} propose to automatically learn a domain decomposition of a PINN by incorporating conditional computing into the PINN framework in order to learn an arbitrary decomposition of the neural network, which is adaptively tuned during the training process. 
Here, conditional computing means that only certain neurons or units of a neural network are activated depending on the network input; see also~\cite{Bengio:2015:CondComputing,Shazeer:2017:ExpertsNN} for more details. Based on the idea of conditional computing the authors introduce GatedPINNs which rely on a number of \textit{experts} which decompose the neural network and are each modeled by a simple neural network themselves. The outputs of the different experts are then aggregated by a linear or nonlinear \textit{gate network}. The authors provide comparative experimental results for the GatedPINN and standard PINN approach and are able to reduce the training time significantly.

A related approach was recently introduced in~\cite{hu2023augmented}, where also a \textit{gate network} is trained to find an augmented domain decomposition and flexibel parameter sharing. Concretely, the gate network does not need any specific interface losses and as output data, provides a weight average of several local subnets. Additionally, it uses a flexible partial parameter sharing of the subnets and initializes the subnets with the hard decomposition obtained from XPINN, resulting in an augmented PINN (APINN).
The authors of~\cite{hu2023augmented} provide results for different PDEs comparing APINN to XPINN and classic PINNs showing that APINNs result in an improved performance and better generalization properties.

In~\cite{dwivedi_distributed_2019}, the authors use ideas from finite volume methods where the computational domain is partitioned into multiple cells to define a distributed PINN (DPINN) for the data-efficient solution of partial differential equations. 
Here, the authors decompose the computational domain into nonoverlapping cells, which can also be interpreted as the subdomains in a DDM, and train a separate, local PINN for each cell with a loss function based exclusively on collocation points from the interior of the cell. Each of these local PINNs is relatively simple and, for the results presented in~\cite{dwivedi_distributed_2019}, only two layers deep. 
In addition to the PINN losses for each cell, a loss term for the interface conditions is introduced which is associated with collocation points located on the boundary of the cells. This interface loss is inspired by the flux conditions of the finite volume method. Then, the DPINN approach is based on training all local PINNs together by minimizing the sum of all local losses plus the interface loss. In particular, the described approach requires the exchange of interface information in each step of the optimization approach used in the training of the network.

A related approach which mostly differs by the treatment of the boundary conditions between neighboring subdomains is presented in~\cite{basir_generalized_2023}. There, the authors use an optimized Schwarz-type nonoverlapping domain decomposition for solving forward and inverse PDE problems using PINNs. Within each subdomain, physics and equality constrained artificial neural networks (PECANNs)~\cite{basir2022physics} are trained in parallel as solvers for each subdomain. Moreover, generalized Robin-type interface transmission conditions are considered as an additional constraint on the solution of each subdomain to ensure continuity. Hence, local information are exchanged at the end of each local training within each outer iteration of the described method. The authors of~\cite{basir_generalized_2023} show that by incorporating the Robin-type interface constraints into the iterative learning process, the effectiveness and reliability of the method can be enhanced.

In~\cite{feeney_breaking_2023}, an end-to-end parallelization of the physics-informed neural PDE solver Mosaic Flow~\cite{wang2022mosaic} based on domain decomposition is presented. 
Roughly speaking, Mosaic Flow is an iterative algorithm inspired by the alternating Schwarz method designed for solving PDEs with arbitrary boundary conditions on diverse domains. It basically consists of two steps. First, local physics-informed PDE solvers (SDNets) are pretrained on small domains with arbitrary boundary conditions, using the discretized boundary function as input data; see also~\cite{wang2022mosaic} for more mathematical details. Then, in a second step, the Mosaic Flow predictor is trained which decomposes a large domain into overlapping subdomains and uses the pretrained SDNets as subdomain solvers.  
In~\cite{feeney_breaking_2023}, the authors combine data parallel training and domain parallelism for inference to develop a distributed domain decomposition algorithm for Mosaic Flow to achieve strong scalability for its parallel training on GPUs. Numerical results for the parallel solution of the Laplace equation on large domains, distributed on 32 GPUs, are presented.

A related approach which is also inspired by the alternating Schwarz method is presented in~\cite{snyder2023domain}. In~\cite{snyder2023domain}, the authors use the alternating Schwarz method and an overlapping domain decomposition to couple local PINNs and full order models, that is, local solutions obtained via the finite element, finite difference, or finite volume method with each other to accelerate the training of PINNs.
Furthermore, three different approaches to enforce Dirichlet boundary conditions within the local subdomain PINNs are empirically compared. Concretely, the three different approaches are the following: A weak boundary condition enforcement through the loss function, a strong enforcement through a solution transformation, and a combination of both such that the boundary conditions are enforced strongly on the computational domain's boundary and weakly on the subdomains' boundaries.
Numerical examples are provided for the one-dimensional steady state advection-diffusion equation with focus on the advection-dominated regime (high Péclet number). The presented results show that the convergence of the Schwarz method is strongly linked to the chosen boundary condition implementation within the local, coupled PINNs.

A different iterative approach to parallelize the training of PINNs while also reducing the communication cost is presented in~\cite{yang_iterative_2023}. In~\cite{yang_iterative_2023}, the authors propose an iteration method which is based on the classical additive Schwarz method. Hence, similarly as in~\cite{Li:2020:DeepDDM,Li:2019:D3M}, an outer training loop for the global network is defined such that, in each iteration of the outer training loop, decoupled smaller local problems are solved which are then stitched back together by a weighted average. In each outer iteration, the local problems obtain their boundary values from the current global iteration and an overlapping decomposition of the collocation points into subdomains is used.
 To compute the next iterate of the global problem, a weighted sum of the previous iterate and the partitioned neural network solution using a relaxation parameter is computed. 
Additionally, to further improve the convergence properties of the described approach, the authors of~\cite{yang_iterative_2023} propose a two-level extension by solving an additional coarse problem with specific structure in each outer iteration.  
Note that related preliminary work has also been published in~\cite{Kim_domain_2022,yang_additive_nodate} by the authors. In~\cite{Kim_domain_2022}, first experimental results for the proposed additive Schwarz algorithm have been presented whereas~\cite{yang_additive_nodate} contains first results regarding the convergence analysis of the approach.

In~\cite{kharazmi_hp-vpinns_2021}, the authors introduce $hp$-variational PINNs ($hp$-VPINNs) which are based on the variational formulation of the residuals of the considered PDE, similar to the approach considered in~\cite{Li:2019:D3M}. In particular, the authors use piecewise polynomial test functions within the variational formulation of the considered PDE such that each test function $v_j$ is always a polynomial of a chosen order over the subdomain $j$ and zero otherwise. 
Let us note that with respect to the implementation of this approach, the authors employ a single DNN to approximate the solution over the whole computational domain despite virtually decomposing the domain into several subdomains; see also \cite[Remark 2.1]{kharazmi_hp-vpinns_2021}. Due to the use of a single DNN no subdomain interface conditions have to be considered in the loss function. However, this makes it more difficult to parallelize $hp$-VPINNs. Different numerical examples of function approximation for continuous and discontinuous functions as well as for solving the Poisson equation are presented.

An approach which is closely related to~\cite{kharazmi_hp-vpinns_2021} and which is also based on the weak, that is, the variational form of a given system of PDEs is presented in~\cite{xu_weak_2021}. The idea of the presented neural networks is to integrate the weak formulation of the PDE in the loss function as well as data constraints and initial or boundary conditions. Due to the theory guided structure of the networks based on the weak formulation of the considered PDE the authors refer to their approach by the acronym TgNN-wf. 
In particular, similar to~\cite{kharazmi_hp-vpinns_2021}, the authors use locally defined test functions to perform a domain decomposition of the computational domain to accelerate the training process. As a further novelty, the authors formulate the original loss minimization problem into a Lagrangian duality problem in order to optimize the penalty term parameters of the different components of the loss function within the training process.
The authors provide comparative results for the strong form TgNN and TgNN-wf for an unsteady-state 2D single-phase flow problem and a 1D two-phase flow problem in terms of accuracy, training time, and robustness with respect to noise in the training data.

In~\cite{sheng_pfnn-2_2022}, the authors use an overlapping domain decomposition approach for the training process of penalty-free neural networks (PFNN), originally introduced in~\cite{sheng2021pfnn}. Applying the ideas introduced in~\cite{sheng2021pfnn} for PFNN, the resulting new PFNN-2 adopts two neural networks, with one network learning the initial condition and essential boundary conditions of a given PDE problem, and a second network approximating the solution of the PDE in its variational form on the remaining part of the domain. To speedup the learning process, the authors introduce an overlapping decomposition of the domain into subdomains and train two local networks, that is, one for the initial and boundary conditions and a second one for the local interior solution, for each subdomain. This approach requires direct exchange of information between neighboring subdomains to enforce continuity via the locally trained networks which learn the boundary conditions of each subdomain.

In~\cite{moseley_finite_2023}, the authors transfer ideas from classical finite element methods where the solution of a PDE or an ODE is expressed as a sum of a finite set of basis functions with a compact support to PINNs. The resulting method uses local window functions to form a global function, that is, the global solution on the entire domain, and is referred to as Finite Basis  Physics-Informed Neural Network (FBPINN). The global function is hence expressed as a sum of the localized window functions which are learned by local neural networks. In particular, the local window or basis functions are defined over small, overlapping subdomains such that, in FBPINNs, no additional interface condition is necessary. In~\cite{moseley_finite_2023}, the authors present numerical results for both small and large, multi-scale problems.

The idea of FBPINNs is further developed in~\cite{dolean_finite_2023} and~\cite{dolean_multilevel_2023}. In~\cite{dolean_finite_2023}, additive, multiplicative, and hybrid iteration methods based on a Schwarz-like domain decomposition method for the training of FBPINNs are introduced. 
In~\cite{dolean_multilevel_2023}, these ideas are extended even further by adding multiple levels of domain decompositions to the solution ansatz, similar as in classical multilevel Schwarz methods. This results in improved convergence properties and increased accuracies in the considered PDE solutions and, at the same time, mitigates convergence problems related to the spectral bias of neural networks.

The performance and convergence properties of L-BFGS~\cite{liu1989limited} optimizers used for the training of PINNs are enhanced in~\cite{kopanicakova_enhancing_2023} by using nonlinear preconditioning strategies based on Schwarz-like domain decomposition approaches. More specifically, nonlinear additive and multiplicative preconditioners are constructed by decomposing the network model, that is, the corresponding parameters in a layer-wise manner. This is, in principle, different from the aforementioned approaches where ideas from DDMs are used to decompose the computational domain. The authors present results for both additive and multiplicative nonlinear preconditioners for different benchmark problems, indicating that the proposed layer-wise decomposition of the network can significantly accelerate the convergence of the L-BFGS optimizer.

\subsection{Domain decomposition for classical supervised and unsupervised machine learning algorithms}
\label{sec:dd_in_others}

\begin{figure}[t]
\centering
\scalebox{0.76}{
\begin{tikzpicture}[line width=1pt]

	\draw[black, rounded corners] (1,1.5) rectangle (4,3.5);
	\draw (2.5,2.5) node[text=black,align=center] {DD for\\classic ML;\\\cref{sec:dd_in_others}};
	
	\draw[thick,->,black] (1,2.5) -- (-1,2.5) -- (-1,1);
	\draw[thick,->,black] (4,2.5) -- (6,2.5) -- (6,1);
	
	\draw[darkblue, rounded corners] (0.5,1) rectangle (-2.5,-1);
	\draw (-1,0) node[text=black,align=center] {supervised ML};
	
	\draw[thick,->,black] (-2.5,0) -- (-4,0) -- (-4,-2);
	\draw[thick,->,black] (0.5,0) -- (2,0) -- (2,-2);
	
	\draw[amber, rounded corners] (4.5,1) rectangle (7.5,-1);
	\draw (6,0) node[text=black,align=center] {unsupervised ML};

	\draw[darkblue, rounded corners] (-5.5,-2) rectangle (-2.5,-3);
	\draw (-4,-2.5) node[text=black,align=center] {neural networks};
	
	\draw[thick,->,black] (-5.5,-2.5) -- (-6.5,-2.5) -- (-6.5,-5.5); 
	\draw[thick,->,black] (-5,-3) -- (-5,-4); 
	\draw[thick,->,black] (-3.5,-3) -- (-3.5,-5.5); 
	\draw[thick,->,black] (-2.5,-2.5) -- (-2,-2.5) -- (-2,-4); 
	
	\draw[darkblue, rounded corners] (0.5,-2) rectangle (3.5,-3);
	\draw (2,-2.5) node[text=black,align=center] {other};
	
	\draw[thick,->,black] (1,-3) -- (1,-4); 
	\draw[thick,->,black] (3,-3) -- (3,-4); 
	
	
	\draw[darkblue, rounded corners] (-7.5,-5.5) rectangle (-5.5,-6.5);
	\draw (-6.5,-6) node[text=black,align=center] {ResNets;\\\cite{gunther_layer-parallel_2019}};
	
	\draw[darkblue, rounded corners] (-6,-4) rectangle (-4,-5);
	\draw (-5,-4.5) node[text=black,align=center] {Fourier NN;\\\cite{malek2023solving}};
	
	\draw[darkblue, rounded corners] (-4.5,-5.5) rectangle (-2.5,-6.5);
	\draw (-3.5,-6) node[text=black,align=center] {RBFNs;\\\cite{mai-duy_mesh-free_2002}};
	
	\draw[darkblue, rounded corners] (-3,-4) rectangle (-1,-5);
	\draw (-2,-4.5) node[text=black,align=center] {CNNs;\\\cite{gu_decomposition_2022,KLW_2023_CNN-DNN}};
	
	\draw[darkblue, rounded corners] (-0.5,-4) rectangle (1.5,-5);
	\draw (0.5,-4.5) node[text=black,align=center] {ELMs;\\\cite{dong_local_2021}};
	
	\draw[darkblue, rounded corners] (2,-4) rectangle (4,-5);
	\draw (3,-4.5) node[text=black,align=center] {GP;\\\cite{park_domain_nodate}};
	
	\draw[thick,->,black] (6,-1) -- (6,-4); 
	
	\draw[amber, rounded corners] (5,-4) rectangle (7,-5);
	\draw (6,-4.5) node[text=black,align=center] {PCA;\\\cite{li_summation_2021}};

\end{tikzpicture}
}
\caption{\label{fig:dd_in_ml} Schematic representation of the articles reviewed in~\cref{sec:dd_in_others} with respect to domain decomposition methods in machine learning.}
\end{figure}
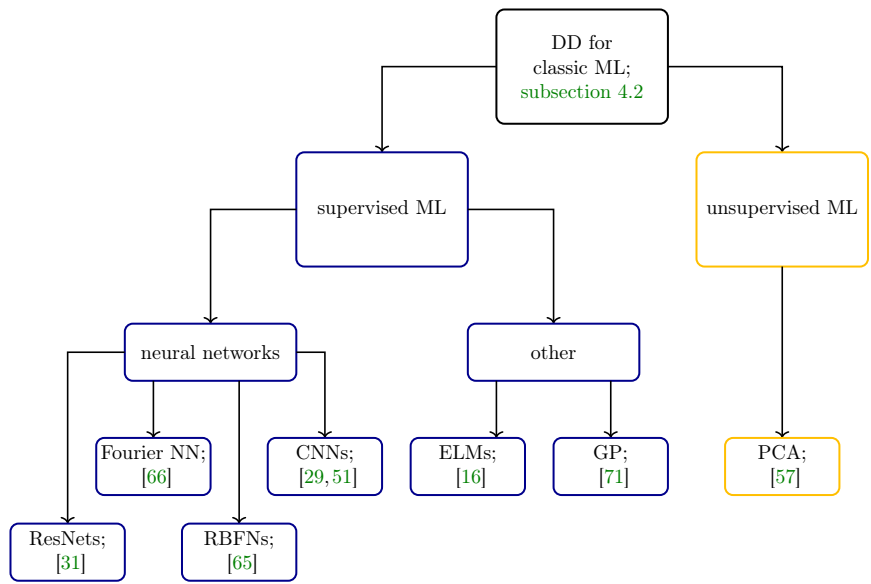

Apart from the wide application area of DDMs for PINNs, the idea of domain decomposition has also been transferred to other machine learning algorithms of which we try to give a broad overview in this section. 
Note that the idea of DDMs is, on the one hand, applied to the training of different supervised machine learning models, for example, convolutional or residual neural networks, and, on the other hand, also to unsupervised algorithms, as, for example, principle component analysis (PCA); see~\cref{fig:dd_in_ml} for a schematic overview. 

In~\cite{gu_decomposition_2022,cai_decomp_2022}, the idea of DDMs is transferred to the training of CNNs for different image recognition problems. The authors propose to decompose a global CNN along its width into a finite number of smaller, local subnetworks which can be trained independently from each other. In particular, the local CNNs are trained in parallel without any communication between the different processes. The obtained weights of the local subnetworks are then used as an initialization for the subsequently trained global network by using a transfer learning strategy. The approach can be interpreted as an analogy to nonlinear preconditioning~\cite{cai2019nonlinear} and the decomposition of the inputs of the global CNN along its width corresponds to a decomposition of the input images into smaller partial images. The authors present results for different two-dimensional image data indicating that the suggested idea can significantly accelerate the training of the global network on a GPU cluster.

A slightly similar but different parallel training approach also for CNNs in image recognition has been proposed in~\cite{KLW_2023_CNN-DNN}. In~\cite{KLW_2023_CNN-DNN}, the authors decompose two- or three-dimensional image data into nonoverlapping or overlapping subsets and train smaller CNNs operating exclusively on parts of the decomposed input data. Based on the decomposed image data, all local CNNs provide a probability distribution with respect to the different classes of the underlying image classification problem. Then, in a second step, a global coarse DNN is trained which evaluates the local decisions into a final, global decision. In particular, this approach does not require the training of a global, large CNN operating on the entire images and the smaller, local CNNs can be trained in parallel on different GPU clusters. The results provided for two- and three-dimensional image data show a significant acceleration of the training and additionally, provide at least the same classification accuracies compared to the training of a single, global CNN for almost all tested datasets.

Another model parallel training approach based on domain decomposition, that is, a layer-parallel training approach has been suggested in~\cite{gunther_layer-parallel_2019} for ResNets~\cite{he2016deep}.
In this work, the authors make use of the fact that the forward propagation through a ResNet can be interpreted as a forward Euler discretization of a time-dependent ordinary differential equation (ODE)~\cite{weinan2017proposal} where the time-dependent control variables represent the weights of the neural network. Hence, the time domain in the time-dependent ODE corresponds to the layer domain of a ResNet. Based on this interpretation, the authors employ a multigrid reduction in time approach~\cite{falgout2014parallel} to replace the classic forward and backward propagation in the training progress by a parallel nonlinear multigrid iteration applied to the layer domain such that chunks of layers can be processed in parallel. In~\cite{gunther_layer-parallel_2019}, scalability and speedup results for this method are presented for three different image recognition problems.

In~\cite{malek2023solving}, domain decomposition techniques are used to define a parallel training process of trigonometric neural networks which use trigonometric activation functions. Based on a nonoverlapping domain decomposition, three local neural networks are trained for each subdomain. The first neural network, denoted as primary neural network is completely local, the second one, denoted as boundary neural network, learns for each subdomain the interface conditions with its neighbors, and the third one, denoted as modifier neural network, then learns the update for the primary network including the interface conditions. Finally, the local solution is defined as the sum of the primary and modifier network. This process is repeated for each subdomain until convergence.

 Another work on the combination of DD and a special type of neural networks is presented in~\cite{mai-duy_mesh-free_2002}. The authors combine different types of radial basis function networks (RBFNs)~\cite{schwenker2001three}, which typically use only one hidden layer,  with the concept of domain decomposition to approximate nonlinear functions or to solve Poisson's equation on a rectangular domain. The computational domain is decomposed into nonoverlapping subdomains and each of the subdomains is discretized by a shallow RBFN. Based on this decomposition, the authors propose an iterative solution algorithm where, in each iteration, all local subproblems are solved using the RBFNs. Then, in each iteration of the outer training loop, the interface condition between the subdomains is estimated and updated using boundary integral equations.

In~\cite{dong_local_2021}, the idea of extreme learning machines (ELMs)~\cite{huang2006extreme} is combined with domain decomposition and local neural networks.  Basically, in order to compute the solution for a given PDE, the authors of~\cite{dong_local_2021} divide the domain into subdomains and represent the solution on each of these subdomains by a local feedforward neural network. To obtain a global solution, continuity constraints are imposed along the subdomain boundaries. Originating from the idea of ELMs, each local neural network consists of a small number of hidden layers, while its last hidden layer is typically wide. As in ELMs, the weights and bias coefficients in all hidden layers of the local neural networks are pre-set to random values and, during training, only the weight coefficients of the output layers of the local neural networks are training parameters. Moreover, also in analogy to ELMs, the overall neural network is trained by a linear or nonlinear least squares computation instead of using the backpropagation algorithm. The authors compare the presented method with the deep Galerkin method (DGM), the PINN method, and classical FEM.

In~\cite{park_domain_nodate}, the idea of domain decomposition is not transferred to a special type of neural network but, instead, to Gaussian Process (GP) regression~\cite{williams2006gaussian}. Given that the computational complexity $O(N^3)$, where $N$ is the number of training samples, is often a major limitation in GP regression, the main idea of the authors is to decompose the optimization problem related to the GP regression into smaller local optimization problems that provide local predictions. This is in analogy of decomposing a domain into smaller subdomains. To deal with the possible mismatch problem of the local predictions along the boundaries of neighboring subdomains, the authors of~\cite{park_domain_nodate} impose continuity constraints, analogously to interface conditions in DDMs.

In~\cite{li_summation_2021}, the authors propose a local PCA approach which also relies on the idea of DDMs. The PCA is applied to different images which are decomposed into nonoverlapping smaller subimages. Here, instead of computing the principal components for the entire images, only principal components for the different subimages are determined. Then, in a second step, the local principal components are used to compute low-dimensional projections of the different subimages which are composed again to a lower dimensional representation of the original global image. The authors show empirically that the decomposed PCA requires a smaller number of principal components to provide an accurate low-dimensional representation of the original image compared to the global PCA. Moreover, the authors provide a mathematical proof that the local PCA reduces the summation pollution effect of the global PCA and thus, results in a lower or equal reconstruction error compared to the global PCA. In other words, the local PCA is never worse than the global PCA; see~\cite{li_summation_2021} for more details.

\section{Machine learning-enhanced domain decomposition methods}
\label{sec:ml_in_dd}

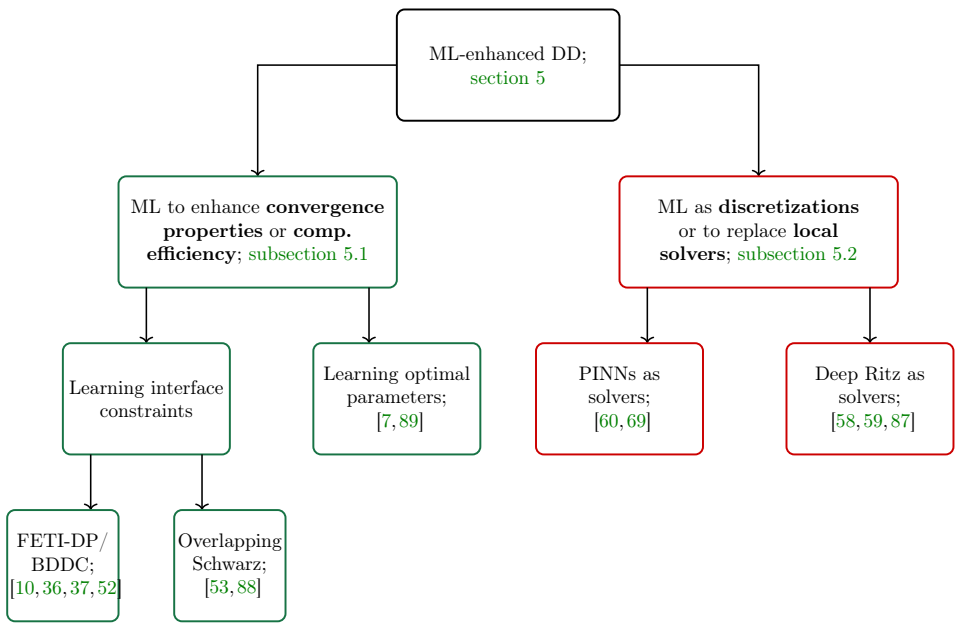
\begin{figure}[t]
\centering
\scalebox{0.74}{
\begin{tikzpicture}[line width=1pt]

	\draw[black, rounded corners] (1,2) rectangle (5,4);
	\draw (3,3) node[text=black,align=center] {ML-enhanced DD;\\\cref{sec:ml_in_dd}};
	
	\draw[thick,->,black] (1,3) -- (-1.5,3) -- (-1.5,1);
	\draw[thick,->,black] (5,3) -- (7.5,3) -- (7.5,1);
	
	\draw[darkspringgreen, rounded corners] (1,1) rectangle (-4,-1);
	\draw (-1.5,0) node[text=black,align=center] {ML to enhance \textbf{convergence}\\\textbf{properties} or \textbf{comp.}\\\textbf{efficiency};~\cref{sec:ml_in_dd_interface}};
	
	\draw[bostonuniversityred, rounded corners] (5,1) rectangle (10,-1);
	\draw (7.5,0) node[text=black,align=center] {ML as \textbf{discretizations}\\or to replace \textbf{local}\\\textbf{solvers};~\cref{sec:ml_in_dd_solvers}};
	
	\draw[thick,->,black] (-3.5,-1)  -- (-3.5,-2);
	\draw[thick,->,black] (0.5,-1)  -- (0.5,-2);
	
	\draw[darkspringgreen, rounded corners] (-2,-2) rectangle (-5,-4);
	\draw (-3.5,-3) node[text=black,align=center] {Learning interface\\constraints};
	
	\draw[darkspringgreen, rounded corners] (-0.5,-2) rectangle (2.5,-4);
	\draw (1,-3) node[text=black,align=center] {Learning optimal\\parameters;\\\cite{taghibakhshi_mg-gnn_2023,frochte}};
	
	\draw[thick,->,black] (-4.5,-4)  -- (-4.5,-5);
	\draw[thick,->,black] (-2.5,-4)  -- (-2.5,-5);
	
	\draw[darkspringgreen, rounded corners] (-4,-5) rectangle (-6,-7);
	\draw (-5,-6) node[text=black,align=center] {FETI-DP/\\BDDC;\\\cite{kimchung,HKLW:2018:ML_acc,HKLW:2020:3D,KLW:2024:JCP}};
	
	\draw[darkspringgreen, rounded corners] (-3,-5) rectangle (-1,-7);
	\draw (-2,-6) node[text=black,align=center] {Overlapping\\ Schwarz;\\\cite{knoke_domain_2023,taghibakhshi_learning_2022}};
	
	\draw[thick,->,black] (5.5,-1)  -- (5.5,-2);
	\draw[thick,->,black] (9.5,-1)  -- (9.5,-2);
	
	\draw[bostonuniversityred, rounded corners] (3.5,-2) rectangle (6.5,-4);
	\draw (5,-3) node[text=black,align=center] {PINNs as\\solvers;\\\cite{Li:2020:DeepDDM,mercier_coarse_2021}};
	
	\draw[bostonuniversityred, rounded corners] (8,-2) rectangle (11,-4);
	\draw (9.5,-3) node[text=black,align=center] {Deep Ritz as\\solvers;\\\cite{Li:2019:D3M,li_deep_2022,sun2022domain}};

\end{tikzpicture}
}
\caption{\label{fig:ml_in_dd} Schematic representation of the articles reviewed in~\cref{sec:ml_in_dd} with respect to machine learning-enhanced domain decomposition methods.}
\end{figure}

In this section, we review recent work which uses different, mainly supervised, machine learning algorithms to enhance different DDMs. In general, one can divide the approaches using ML to enhance DDMs into two classes; see also~\cite{HKLW:2020:GammReview}. 
The first class consists of methods where ML is used to improve the convergence properties or the computational efficiency within classical DDMs. Typically, this is done by learning or approximating optimal interface conditions or by learning other optimal parameters. Existing approaches of this class are considered in~\cref{sec:ml_in_dd_interface}.
In the second class, different types of neural networks are used as discretization methods and replace classical local subdomain or coarse solvers based on finite elements or finite differences. Existing methods of this class are summarized in~\cref{sec:ml_in_dd_solvers}. For a schematic overview of the existing literature, see~\cref{fig:ml_in_dd}.

\subsection{Machine learning to enhance convergence properties or computational efficiency of domain decomposition methods}
\label{sec:ml_in_dd_interface}

As explained in~\cref{sec:dd}, introducing appropriate global communication and interface conditions 
between subdomains in a DDM is a crucial task to ensure the correct and globally continuous solution as well as a satisfactory convergence behavior. 
Hence, a number of different attempts have been made to automatically construct such interface conditions, for example, by using regression neural networks.

In~\cite{kimchung}, the authors propose an approach to learn an adaptive coarse space for BDDC (Balancing Domain Decomposition by Constraints) algorithms for stochastic elliptic PDEs. The authors focus on a specific adaptive coarse space~\cite{kim1,kim2017bddc} where local generalized eigenvalue problems have to be solved for parts of the interface between neighboring subdomains. In particular, these eigenvalue problems depend on the (stochastic) coefficients of the PDE. Hence, the authors train a dense feedforward neural network to approximate the relation between the stochastic PDE coefficients and the adaptive coarse space. As an input for the neural networks, a truncation of the Karhunen-Lo\`{e}ve expansion~\cite{zhang2004efficient} is used. Numerical results are presented for different oscillatory and high contrast coefficients, a fixed discretization of the domain $\Omega \subset \mathbb{R}^2$, and a neural network with a single hidden layer. 
Note that, so far, this method is not independent of the discretization and thus has to be retrained for different mesh resolutions.

As a preparatory step to automatically construct an adaptive coarse space for DDMs, in~\cite{HKLW:2018:ML_acc,HKLW:2020:3D}, the authors present a supervised classification approach to identify parts of the interface where adaptive coarse constraints are necessary. Then, adaptive constraints resulting from the local generalized eigenvalue problem in~\cite{sousedik1,sousedik2} are integrated into a FETI-DP (Finite Element Tearing and Interconnecting - Dual Primal) coarse space exclusively for edges or faces, respectively, which are categorized as critical by a pretrained neural network. In particular, the decision of the neural network is based on a sampling procedure of the coefficient function of the underlying elliptic PDE and is independent of the finite element discretization. Hence, the trained neural network can be evaluated for different finite element meshes as used in the training procedure and can be applied to both, two- and three-dimensional problems.
Additionally, in~\cite{HKLW:2020:GDSW}, the described approach is also analogously applied to the AGDSW (Adaptive Generalized Dryja-Smith-Widlund) method~\cite{Heinlein:2018:AGD} in two spatial dimensions.

In~\cite{KLW:2024:JCP}, the work described above has been extended to directly learn the adaptive constraints to setup a robust FETI-DP coarse space in two spatial dimensions. Here, the authors train a fixed number of regression neural networks to directly predict an approximation of adaptive edge constraints. Again, the input data used for the neural networks rely on a mesh-independent sampling strategy such that the trained neural networks can be applied to different finite element mesh resolutions. Numerical results are provided for different elliptic PDE problems and realistic coefficient distributions orginating from a microsection of a dual-phase steel. In~\cite{KLW:2022:DD27_TR}, these results are further extended to irregular decompositions obtained by the graph partitioning software METIS~\cite{METIS}.

In order to provide a robust and efficient iterative solution method for the time-harmonic Maxwell's equations, in~\cite{knoke_domain_2023}, the authors use a feedforward neural network to approximate the interface operator in an optimized Schwarz method~\cite{dolean2009optimized}. The concrete formulation of the optimized Schwarz method for the time-harmonic Maxwell's equations is based on~\cite{el2015quasi}. In~\cite{knoke_domain_2023}, preliminary numerical results for a two-dimensional domain decomposed into two subdomains are provided. Additionally, comparative results for different network architectures and different activation functions are presented and compared with classical DDMs.

A related approach to learn optimal interface conditions also for an optimized Schwarz method~\cite{gander2001optimized} using Graph Convolutional Neural Networks (GNNs) is suggested in~\cite{taghibakhshi_learning_2022}. 
 There, the authors aim to learn the subdomain interface matrices for both, structured grids where the optimal matrix values can be predicted by Fourier analysis and for unstructured grids where the optimal Schwarz parameters are generally not known.
Hence, the GNNs are trained using as input the degrees of freedom of a discretized problem, its decomposition, and sparsity constraints on the subdomain interface matrices. As output, the networks predict approximated optimized values for these matrices. Additionally, the authors consider an improved loss function which is a relaxation of the ideal Frobenius norm minimization.

In~\cite{taghibakhshi_mg-gnn_2023}, a new multigrid GNN (MG-GNN) architecture to learn optimized parameters in two-level DDMs is introduced. In particular, the new MG-GNN model is used to learn the Robin-type subdomain boundary conditions in a two-level optimized Schwarz method and the overall coarse-to-fine interpolator. In a way, the work in~\cite{taghibakhshi_mg-gnn_2023} can thus be seen as an extension of the method in~\cite{taghibakhshi_learning_2022}, where also GNNs are used to learn optimal interface conditions in a one-level optimized Schwarz method.

A slightly different approach, where not adaptive constraints themselves but interface conditions in form of varying overlap between neighboring subdomains are learned is presented in~\cite{frochte}. 
More precisely, a method is suggested to optimize the width $\delta$ of the overlap in Schwarz methods such that the number of expected floating point operations until convergence is minimized.  
The authors compare four different regression algorithms for two-dimensional diffusion problems with jumps in the coefficient function. 
 As input for the different regression algorithms, the authors use certain features collected for each subdomain as well as its neighborhood, as, for example, the maximal or minimal coefficient within certain sets of rows and columns of degrees of freedom in the surrounding of the boundary of the overlapping subdomain.

\subsection{Replacing the subdomain solvers in domain decomposition methods by neural networks}
\label{sec:ml_in_dd_solvers}

 In addition to learning optimal interface conditions to improve the robustness of DDMs, a different approach that has attracted much attention very recently is to replace the subdomain solvers by machine learning algorithms as, for example, different types of neural networks. 
As a starting point for this section, we consider again boundary value problem~\eqref{eq:bvp} and, in the following, review some works which aim to solve the local subdomain problems by training a neural network $\mathcal{N}$ with certain collocation points and the loss function given in~\eqref{eq:loss} or~\eqref{eq:loss_vp3}, respectively.

In~\cite{Li:2020:DeepDDM} and~\cite{Li:2019:D3M}, the subdomain solvers in a parallel overlapping Schwarz method are replaced by PINNs and the Deep Ritz method, respectively, such that each of these local networks solves a single subdomain problem. The latter method is named D3M in~\cite{Li:2019:D3M} and the former method DeepDDM in~\cite{Li:2020:DeepDDM}. 
As explained above, in any DDM, it is crucial to exchange information between the local subdomains in order to obtain a continuous, global solution. 
In DeepDDM~\cite{Li:2020:DeepDDM} as well as in D3M~\cite{Li:2019:D3M}, which are both based on a parallel overlapping Schwarz fixed point iteration, the exchange of information is enforced by additional boundary conditions, which change in each  fixed-point iteration until convergence; see~\cite{Lions:1988:SAM} or~\cite{toselli} for further details on the parallel overlapping Schwarz method.
These additional boundary conditions result in a third loss term of the local subdomain networks with additional collocation points for the exchange of information between neighboring subdomains.

In~\cite{mercier_coarse_2021}, the authors extend the DeepDDM algorithm~\cite{Li:2020:DeepDDM} described above by a coarse space correction, similarly to what is done in traditional two-level DDMs~\cite{toselli}.
The aim is to improve the convergence properties of the resulting algorithm for an increasing number of subdomains due to the enhanced exchange of global information. To define a coarse problem, the authors of~\cite{mercier_coarse_2021} introduce an additional coarse sampling of the complete domain defined by a set of coarse collocation points. These coarse collocation points are then used as input data for an additional fully connected network that is trained subsequently after the local networks in each outer training iteration. 
The authors present numerical results for different elliptic boundary value problems indicating the advanced numerical scalability for an increasing number of subdomains compared to the original DeepDDM approach. 

A different extension of the D3M approach~\cite{Li:2019:D3M} is presented in~\cite{li_deep_2022}. Here, the aim is to mitigate the spectral bias~\cite{rahaman2019spectral,luo2019theory,ronen2019convergence} observed with the local networks in the original D3M approach by introducing new multi Fourier feature networks (MFFNet) in each local subdomain. This results in a new Fourier feature based deep DDM (F-D3M) which is more successful in approximating the high frequency modes of a PDE solution. In particular, in order to obtain satisfactory approximations for both, low and high frequency modes, the MFFNet consists of different subnetworks with different random frequencies of the Fourier features and whose outputs are concatenated with a linear layer to compute the final output of the MFFNet.

In~\cite{sun2022domain}, the subdomain solvers within a DDM are also replaced by a Deep Ritz approach but, in contrast to~\cite{Li:2019:D3M}, the computational domain is decomposed into nonverlapping subdomains. In contrast to replacing the subdomain solvers by neural networks within an overlapping Schwarz method as done in D3M or DeepDDM, in its nonoverlapping counterpart a direct flux exchange across subdomain interfaces is required. In particular, simply replacing the subdomain solvers in nonoverlapping subdomains by the classic Deep Ritz approach often results in a local minimizer where the given Dirichlet boundary conditions are satisfied satisfactory but with inaccurate Neumann boundary conditions; see~\cite{sun2022domain}. In order to address the issue of the low accuracy of the related flux transmission between neighboring nonoverlapping subdomains, the authors of~\cite{sun2022domain} propose a compensated Deep Ritz method as subdomain solvers that enable reliable flux transmission in the presence of erroneous interface conditions. The compensated Deep Ritz method defines a weighted sum of the Ritz form of the loss function and the PINN residual loss given that PINNs are empirically found to provide more accurate estimates of partial derivatives. The authors provide numerical examples for elliptic problems and two and four subdomains.

\section{Discussion and future work}
\label{sec:concl}

In this paper, we have explored how established techniques and ideas from DDMs can effectively be combined with machine learning algorithms to develop major advances in SciML methods. In particular, we have categorized all reviewed work in the two main classes of DDMs for machine learning and machine learning-enhanced DDMs. 
By taking all cited work into account, we observe that, so far, higher efforts exist to implement ideas from DDMs into machine learning models than vice versa. However, in both areas, several ongoing challenges and opportunities for further developments and research exist.

In the field of DDMs for machine learning, a very large group of publications exists that use DDMs within physics-aware neural networks. Here, we have observed a strong focus of using DDMs as model or data parallelism to efficiently accelerate and distribute the training of physics-aware neural networks, or, in particular, PINNs. Furthermore, decomposing the training of one large global neural network into training several smaller neural networks in parallel can also have important impacts on the generalization properties of the underlying machine learning model. As investigated analytically in~\cite{hu_when_2022}, decomposing the global optimization problem into smaller subproblems decreases the complexity of each subproblem which can improve the generalization of the resulting global model. On the other hand, each subproblem is trained with a smaller amount of training data which tends to increase overfitting of the local models. Thus, a tradeoff between both effects has to be considered. 
Additionally, training decoupled local PINN models in a preliminary phase which are subsequently aggregated into a global model can also help to mitigate the spectral bias often observed for deep neural networks.

Besides, concepts from DDMs have also been integrated into classical machine learning models, for example, CNNs for image recognition, ResNets, Gaussian Processing, or PCA. Here, the primary objective is also to distribute the training of the respective machine learning model to different GPUs in order to accelerate the required training time. As a secondary effect, also for DDMs in classical machine learning models, in some cases, an improvement of the accuracy of the trained model can be observed. For some models, a mathematical comparison of the original global model versus the decomposed model already exists, for example,~\cite{li_summation_2021,gu_decomposition_2022}, whereas for other cases, a mathematical proof of a corresponding correlation is still an open research question. 

Within the second reviewed field of machine-learning enhanced DDMs, a range of work exists on learning interface conditions or restriction operators by regression neural networks. So far, all of these approaches are developed for concrete domain decomposition methods and test problems, for example, for FETI-DP and BDDC or optimized Schwarz methods for stationary diffusion problems, stochastic PDEs, or Maxwell's equations, respectively. Here, a possible priority for future work could be the aim to develop more universal approaches which can be applied to various problems and DDMs or are independent of the underlying finite element discretization. Currently, only some of the cited approaches partly fulfill such requirements. 
Additionally, some work on learning optimal parameters in DDMs exist, as, for example, learning the optimal width of overlap for overlapping subdomains, optimal stopping criterions, or optimal restriction operators. These efforts, considering the current state, are also restricted to very concrete cases where specific information with respect to the considered DDM and test problem are used as input data for different regression models. Thus, in this area, more generally valid approaches are also a topic of potential future research. 

Finally, there are several efforts to use PINNs or Deep Ritz methods to replace the subdomain solvers or the discretization within classical DDMs. These approaches make use of the universal approximation capabilities of neural networks to represent, under certain assumptions on the activation function, etc., solutions of PDEs. However, there are some challenges within these approaches in order to be competitive with classical numerical methods for the solution of PDE problems. Most notably, most existing approaches require to be retrained for varying boundary and initial conditions which makes the training of such approaches computationally expensive and limits the applicability of the resulting models. Hence, practically relevant applications for such models have to be carefully thought of, for example, large-scale multi-physics problems, and more flexible surrogate models will much likely be an important topic of future research.

Neural PDE solvers or, more general, neural operators can offer a remedy for the described drawback that PINNs or Deep Ritz usually need to be retrained for varying boundary conditions. Neural operators are a type of model that is trained to approximate the PDE solution operator and to be evaluated for various instances of a boundary value problem with different boundary conditions. Practically, this means that the neural operator is trained with a representative set of discretized boundary functions as input data.
With respect to the aim of applying neural operators to large, complex domains with arbitrary boundary conditions, some efforts have been made to combine neural operators with ideas from DDMs to enable the efficient solution of complex boundary value problems and a scalable training of these models. The different avenues of ongoing research, so far, still face several challenges, for example, with respect to applicability to irregular meshes, complex PDE solutions, and computationally expensive training. Hence, combining DDMs and other promising concepts as, for example, transformers with neural operators~\cite{hao2023gnot} to address these challenges is likely to be a future research topic of rising interest and importance.

\bibliographystyle{siamplain} 
\bibliography{Review_Paper} 

\end{document}